\documentclass[11pt, leqno, twoside]{article}
\usepackage{}

\usepackage{amsfonts}

\usepackage{amssymb}
\usepackage{amsmath}
\usepackage{amsthm}
\usepackage{xcolor}
\usepackage{mathrsfs}
\usepackage{txfonts}
\usepackage{extarrows}

\usepackage{indentfirst}

\allowdisplaybreaks

\def\fz{\infty}

\def\dz{\delta}

\def\lz{\lambda}
\def\nn{{\mathbb{N}}}
\def\vz{\varphi}
\def\Oz{\Omega}

\def\lf{\left}
\def\r{\right}

\def\lv{{L^{\vz}(K)}}
\def\mz{{L^{0}(K)}}
\def\cl{{\Oz_{\vz}(K)}}
\def\clp{{\Oz_{\phi}(K)}}
\def\clz{{\Oz_{\vz_1}(K)}}
\def\clzz{{\Oz_{\vz_2}(K)}}
\def\clzzz{{\Oz_{\vz_3}(K)}}
\def\mk{{\mathcal{M}(K)}}

\def\noz{\nonumber}

\def\supp{\mathop\mathrm{supp}}
\def\esup{\mathop\mathrm{\,ess\,sup\,}}

\newtheorem{theorem}{Theorem}[section]
\newtheorem{lemma}[theorem]{Lemma}

\newtheorem{question}[theorem]{Question}
\theoremstyle{definition}
\newtheorem{remark}[theorem]{Remark}
\newtheorem{definition}[theorem]{Definition}
\renewcommand{\appendix}{\par
   \setcounter{section}{0}%
   \setcounter{subsection}{0}%
   \setcounter{subsubsection}{0}%
   \gdef\thesection{\@Alph\c@section}%
   \gdef\thesubsection{\@Alph\c@section.\@arabic\c@subsection}%
   \gdef\theHsection{\@Alph\c@section.}%
   \gdef\theHsubsection{\@Alph\c@section.\@arabic\c@subsection}%
   \csname appendixmore\endcsname
 }

\numberwithin{equation}{section}

\begin{document}

\arraycolsep=1pt

\title{\bf\Large Some Sets of First Category in Product
Calder\'{o}n--Lozanovski\u{\i} Spaces on Hypergroups\footnotetext {\hspace{-0.35cm}
2020 {\it Mathematics Subject Classification}. Primary 46E30;
Secondary 43A62, 43A15, 54E52.
\endgraf {\it Key words and phrases}. Locally compact hypergroup, Calder\'{o}n--Lozanovski\u{\i} space,
convolution, first category, Orlicz space.
\endgraf
This project is supported by the National
Natural Science Foundation of China (Grant Nos.~12371102 and 12371025).
}}
\author{Jun Liu, Yaqian Lu and Chi Zhang\footnote{Corresponding author, E-mail: zclqq32@cumt.edu.cn}}
\date{}
\maketitle

\vspace{-0.9cm}

\begin{center}
\begin{minipage}{13cm}
{\small {\bf Abstract}\quad
Let $K$ be a locally compact hypergroup with a left Haar measure $\mu$ and $\Omega$ be a Banach ideal
of $\mu$-measurable complex-valued functions on $K$. For Young functions $\{\varphi_i\}_{i=1,2,3}$,
let $\Omega_{\varphi_i}(K)$ be the corresponding Calder\'{o}n--Lozanovski\u{\i} space associated with
$\varphi_i$ on $K$. Motivated by the remarkable work of Akbarbaglu et al. in [Adv. Math. 312 (2017), 737-763],
in this article, the authors give several sufficient conditions for the sets
$$\left\{(f,g)\in\Omega_{\varphi_1}(K)\times\Omega_{\varphi_2}(K):\
|f|\ast |g|\in\Omega_{\varphi_3}(K)\right\}$$
and
$$\left\{(f,g)\in\Omega_{\varphi_1}(K)\times\Omega_{\varphi_2}(K):\
\exists\,x\in U,\ (|f|\ast |g|)(x)\ \text{is\ well\ defined}\right\}$$
to be of first category in the sense of Baire, where $U\subset K$ denotes a compact set.
All these results are new even for Orlicz(--Lorentz) spaces on hypergroups.}
\end{minipage}
\end{center}

\section{Introduction}\label{s1}

Let $G$ be a locally compact group with a fixed left Haar measure $\mu$
and $L^p(G)$ be the Lebesgue space on $G$ with $p\in[1,\fz]$.
For two $\mu$-measurable functions $f$ and $g$ on $G$, their convolution is defined as
$$(f\ast g)(x):=\int_{G}f(y)g(y^{-1}x)\,d\mu(y)$$
for every $x\in G$ such that the function $y\mapsto f(y)g(y^{-1}x)$ is $\mu$-integrable.
Moreover, $f\ast g$ is said to exist if $(f\ast g)(x)$ exists for almost every $x\in G$.

Let us give a brief review about the research on $L^p$-conjecture. This famous conjecture was
first formulated by Rajagopalan in his Ph.D. thesis in 1963 as follows.

\setcounter{theorem}{0}
\renewcommand{\thetheorem}{\Alph{theorem}}

\begin{theorem}\label{T1}
Let $G$ be a locally compact group and $p\in(1,\fz)$. Then $f\ast g$ exists and belongs to
$L^p(G)$ for any $f,\,g\in L^p(G)$ if and only if $G$ is compact.
\end{theorem}

To be precise, the first well-known result related to this conjecture was due to \.{Z}elazko \cite{z61}
and Urbanik \cite{u61} in 1961, in which they proved that Theorem \ref{T1} is true for any locally
compact abelian group $G$. Later, \.{Z}elazko \cite{z63} and Rajagopalan \cite{Ra66} independently
gave the proof of Theorem \ref{T1} for $p\in(2,\fz)$; see also Rajagopalan's works \cite{Ra63} for the case
when $p\in[2,\fz)$ and $G$ is discrete, \cite{Ra66} for the case when $p=2$ and $G$
is totally disconnected, and \cite{Ra67} for the case when $p\in(1,\fz)$ and $G$ is either nilpotent
or a semi-direct product of two locally compact abelian groups as well as Rickert's work \cite{Ri68}
for the case when $p=2$ and $G$ is as in Theorem \ref{T1}. In the joint work \cite{rz65},
the aforementioned two authors established the truth of Theorem \ref{T1} for any locally compact solvable
group $G$; see also Greenleaf's monograph \cite{g69}. After thirty years of struggling, until 1990, Saeki
\cite{Sa90} gave an affirmative answer to this conjecture by quite elementary means.

One significant topic on $L^p$-conjecture is the study of its quantitative version, which originated from  G{\l}\c{a}b and Strobin \cite{gs10}. Concretely, they obtained the following conclusion and also raised
several interesting questions; see \cite[Theorem 1 and Section 4]{gs10}.

\begin{theorem}\label{T2}
Assume that $G$ is a non-compact but locally compact group. Let $p,\,q\in(1,\fz)$ satisfy
$1/p+1/q<1$. Then, for each compact set $U\subset G$, the set
$$\left\{(f,g)\in L^p(G)\times L^q(G):\
\exists\,x\in U,\ (f\ast g)(x)\ is\ well\ defined\right\}$$
is $\sigma$-lower porous (see \cite[p.\,584]{gs10} for its precise definition).
\end{theorem}

To answer G{\l}\c{a}b and Strobin's questions, Akbarbaglu and Maghsoudi \cite{am12} confirmed that, if $G$
is a non-unimodular locally compact group, then, for each compact set $U\subset G$, the set
$\{(f,g)\in L^p(G)\times L^q(G):\ \exists\,x\in U,\ (f\ast g)(x)\ \text{is\ well\ defined}\}$
is $\sigma$-lower porous, where $p\in(1,\fz)$ and $q\in[1,\fz)$. The same assertion for the case when $G$
is non-discrete with $p\in(0,1)$ and $q\in(0,\fz]$ was explored by Akbarbaglu and Maghsoudi \cite{am13}.
Moreover, in \cite{am13}, they also showed that, if $G$
is a non-discrete locally compact group, then, for each compact set $U\subset G$, the set
$$\lf\{(f,g)\in L^p(G)\times L^q(G):\ f\ast g\in L^r(U,\mu|_U)\r\}$$
is $\sigma$-lower porous, where $p,\,q\in[1,\fz)$ and $r\in[1,\fz]$ such that $1/p+1/q>1+1/r$. This
meanwhile solved an open problem proposed by Saeki \cite{Sa90}.

Another meaningful line of inquiry into $L^p$-conjecture is extending the space $L^p$ to some more general
function spaces on a locally compact group $G$ such as the Orlicz space $L^{\vz}(G)$ etc.
Precisely, the aim of this
research direction is to find out what needs to assumed on Young functions $\{\varphi_i\}_{i=1,2,3}$
such that $f\ast g\in L^{\vz_3}(G)$ provided $f\in L^{\vz_1}(G)$ and $g\in L^{\vz_2}(G)$, where, for each
$i\in\{1,2,3\}$, $L^{\vz_i}(G)$ denotes the corresponding Orlicz space associated with $\varphi_i$ on $G$.
Recall that O'Neil \cite{Ne65} proved that, if $G$ is a unimodular locally compact group and the Young
functions $\{\varphi_i\}_{i=1,2,3}$ satisfy $\vz_1^{-1}(x)\vz_2^{-1}(x)\le x\vz_3^{-1}(x)$ for any
$x\in[0,\fz)$, then $f\ast g\in L^{\vz_3}(G)$ for any $f\in L^{\vz_1}(G)$ and $g\in L^{\vz_2}(G)$,
that is, the convolution operator is bounded from $L^{\vz_1}(G)\times L^{\vz_2}(G)$ into $L^{\vz_3}(G)$.
In addition, Hudzik, Kami\'{n}ska and Musielak in \cite[Theorem 2]{hkm85} presented some equivalent
conditions for an Orlicz space to be a convolution Banach algebra:

\begin{theorem}\label{T3}
Let $G$ be a locally compact abelian group and $\vz$ a Young function satisfying the so-called
$\Delta_2$-condition. Then $f\ast g$ exists and belongs to $L^{\vz}(G)$ for any $f,\,g\in L^{\vz}(G)$
if and only if $G$ is compact or $\displaystyle\lim_{x\to0}\frac{\vz(x)}{x}>0.$
\end{theorem}

This result was later extended, respectively, by Kami\'{n}ska and Musielak \cite{km89} to the form of
$L^{\vz_1}(G)\ast L^{\vz_2}(G)\subset L^{\vz_3}(G)$ and by Akbarbaglu and Maghsoudi \cite{am13'}
to the amenable group. In particular, in \cite{am13'}, they obtained that, if $G$ is amenable,
$\vz$ is a $\Delta_2$-regular $N$-function and $f\ast g\in L^{\vz}(G)$ for any $f,\,g\in L^{\vz}(G)$,
then $G$ is compact. It now leads us naturally to the question: Under the same assumptions, does
the compactness of $G$ imply that $f\ast g\in L^{\vz}(G)$ for any $f,\,g\in L^{\vz}(G)$?
And more generally, when $f\ast g$ exists as a well-defined function on $G$ for any $f,\,g\in L^{\vz}(G)$?
A discussion on this question can be found in \cite{am12'} of Akbarbaglu and Maghsoudi, which is considered
as the initiation on the research of the quantitative $L^p$-conjecture for Orlicz spaces. The main result
of \cite{am12'} is as follows; see \cite[Theorems 2.1 and 2.4]{am12'}.

\begin{theorem}\label{T4}
Assume that $G$ is a non-unimodular locally compact group, $\vz_1$ is an $N$-function, and $\vz_2$ is a
finite Young function. Then, for each compact set $U\subset G$, the set
$$\left\{(f,g)\in L^{\vz_1}(G)\times L^{\vz_2}(G):\
\exists\,x\in U,\ (f\ast g)(x)\ is\ well\ defined\right\}$$
is $\sigma$-lower porous. The conclusion is also true for the non-compact locally compact group with
an additional condition on $\vz_1$ and $\vz_2$:
\begin{align*}
\liminf_{x\to0}\frac{\vz_1^{-1}(x)\vz_2^{-1}(x)}{x}=\fz.
\end{align*}
\end{theorem}

An improvement and unification of these above results was completed by
Akbarbaglu et al. in \cite{agms} under a more general setting of Calder\'{o}n--Lozanovski\u{\i} spaces,
which include many familiar function spaces as special cases, such as Orlicz spaces, Orlicz--Lorentz spaces
and $p$-convexification spaces etc. Recall that the Calder\'{o}n--Lozanovski\u{\i} space was first
introduced by Calder\'{o}n in \cite{Ca64} and well developed by Lozanovski\u{\i} in \cite{Lo65,Lo73}.
They play a key role in the theory of interpolation since their fine structure is the interpolation
functor for positive operators and, also for all linear operators with some additional assumptions,
like the Fatou property or separability. For more development on this space, we refer the
reader to, for instance, \cite{fk24,gms14,hsy21,kl10,klm14,kmm20,lm17,Ma89,rt18}.
Let $G$ be a locally compact group and, for each $i\in\{1,2,3\}$, let
$\varphi_i$ be a Young function, $E$ an ideal of measurable complex-valued functions defined on $G$
and $E_{\varphi_i}(G)$ be the corresponding Calder\'{o}n--Lozanovski\u{\i} space. To be precise,
the authors of \cite{agms} investigated when the sets
$$\left\{(f,g)\in E_{\varphi_1}(G)\times E_{\varphi_2}(G):\
|f|\ast |g|\in E_{\varphi_3}(G)\right\}$$
and
$$\left\{(f,g)\in E_{\varphi_1}(G)\times E_{\varphi_2}(G):\
\exists\,x\in V,\ (f\ast g)(x)\ \text{is\ well\ defined}\right\}$$
are of first category in the sense of Baire, where $V\subset G$ denotes a compact set.
The analogue quantitative problems for pointwise products in the context of complete $\sigma$-finite
measure spaces have also been considered in \cite{agms}.

On another hand, the discovery of the notion of hypergroups, which has a group-like structure that can
be expressed in terms of an abstract convolution of measures, is one of the decisive steps in studying
stochastic dynamical systems within the framework of algebraic-topological structures specified by
invariance conditions. Moreover, due to the celebrated work of Dunkl \cite{Du73} and Jewett \cite{Je75}
on the fundamental theory of both generalized translation spaces and hypergroups, there has been an
increasing interest in extending harmonic analysis and representation theory from locally compact
semigroups to hypergroups; see, for instance, \cite{bh95} and the recent monograph \cite{La23}.
Particularly, we refer the reader to \cite{th18,bkt22} for the study on $L^p$-conjecture in the setting
of hypergroups and to \cite{tba24} for its quantitative version. Then a natural question is to explore
the results related to the quantitative $L^p$-conjecture for Orlicz spaces on hypergroups,
and more generally, for Calder\'{o}n--Lozanovski\u{\i} spaces on hypergroups.

Motivated by the aforementioned work, in this article, we give several sufficient conditions for
some subsets of product Calder\'{o}n--Lozanovski\u{\i} spaces on hypergroups to be of first category in
the sense of Baire. The main ideal we used to prove our results comes from \cite{agms,tba24}.
The difficulty for the present setting of hypergroups stems from the fact that, compared with locally
compact groups, the convolution of two Dirac measures of a hypergroup is not necessarily a Dirac measure.
To overcome this difficulty, we have to find some mild conditions on hypergroups which are automatically
true for locally compact groups. The remainder of this article is organized as follows.
In Section \ref{s2}, we recall the notions concerning hypergroups, Orlicz spaces,
Calder\'{o}n--Lozanovski\u{\i} spaces etc., and then state the main results of this article.
Sections \ref{s3} and \ref{s4} are devoted to proving these main results.

Hereinafter, let $\nn:=\{1,2,\ldots\}$ and,
for any set $E$, denote by $\mathbf{1}_E$ its \emph{characteristic function}.
For any complex function $f$, denote by ${\rm Re}(f)$ its real part.

\section{Preliminaries and main results}\label{s2}

In this section, we first recall some notions
which will be used throughout this article and then state the main results of this article.

\subsection{Hypergroups}\label{s2.1}

We commence with the definition of hypergroups from \cite{bkt22,th18}, which was initially introduced in
\cite{Du73,Je75}. Let $K$ be a locally compact Hausdorff space. Denote by $\mathcal{M}(K)$ the space
of all (complex-valued) Radon measures on $K$ and by $\mathcal{M}^+(K)$ the set of all
non-negative measures in $\mathcal{M}(K)$. The support of each $\mu\in\mathcal{M}(K)$ is
denoted by $\supp(\mu)$, and for any $x\in K$, the symbol $\dz_x$ denotes the Dirac measure at $x$.

\setcounter{theorem}{0}
\renewcommand{\thetheorem}{\arabic{section}.\arabic{theorem}}

\begin{definition}\label{2d1}
Let $K$ be a locally compact Hausdorff space. If there are a convolution $\ast:\ \mk\times\mk\to\mk$,
an involution $x\mapsto x^-$ on $K$, and an element $e\in K$ (called the identity element) satisfying
the following conditions:
\begin{enumerate}
\item[{\rm(i)}] $(\mk,+,\ast)$ is a complex Banach algebra;
\item[{\rm(ii)}] for any non-negative $\mu,\nu\in\mk$, $\mu\ast\nu$ is also a non-negative measure in
$\mk$ and the mapping $(\mu,\nu)\mapsto\mu\ast\nu$ from $\mk^+\times\mk^+$ into $\mk^+$ is continuous
with respect to the cone topology, that is, the weakest topology such that, for each
$f\in C_{\rm c}^+(K)$, the mappings $\mu\mapsto\int_Kf(x)\,d\mu(x)$ and $\mu\mapsto\mu(K)$
are continuous, where $C_{\rm c}^+(K)$
denotes the set of all the non-negative continuous complex-valued functions on $K$ with compact support;
\item[{\rm(iii)}] for any $x, y\in K$, $\dz_x\ast\dz_y$ is a probability measure with compact support,
that is, $\dz_x\ast\dz_y\in\mk^+$, $(\dz_x\ast\dz_y)(K)=1$ and $\supp(\dz_x\ast\dz_y)$ is compact;
\item[{\rm(iv)}] for any $x\in K$, $\dz_e\ast\dz_x=\dz_x=\dz_x\ast\dz_e$;
\item[{\rm(v)}] the mapping $x\mapsto x^-$ on $K$ is a homeomorphism, and for any $x, y\in K$,
$(x^-)^- =x$ and $(\dz_x\ast\dz_y)^-=\dz_{y^-}\ast\dz_{x^-}$, where, for any $\mu\in\mk$ and
$\Omega\subset K$, $\mu^-(\Omega):=\mu(\{x^-:\ x\in\Omega\})$.
Also, for any $x, y\in K$, $e\in\supp(\dz_x\ast\dz_y)$ if and only if $x=y^-$;
\item[{\rm(vi)}] the mapping $(x,y)\mapsto\supp(\dz_x\ast\dz_y)$ from $K\times K$ into $\mathfrak{C}(K)$
is continuous, where $\mathfrak{C}(K)$ denotes the space of all the non-empty compact subsets of $K$
equipped with the Michael topology whose subbasis has the form
$\{H\in\mathfrak{C}(K):\ H\cap U\neq\emptyset\ {\rm and}\ H\subset V\},$
in which $U$ and $V$ are open subsets of $K$.
\end{enumerate}
Then the quadruple $(K,\ast,^-,e)$ is called a \emph{locally compact hypergroup}
(or simply a \emph{hypergroup}).
\end{definition}

\begin{remark}\label{2r1}
\begin{enumerate}
\item[{\rm(i)}] When there is no confusion, we write $K$ instead of $(K,\ast,^-,e)$.
\item[{\rm(ii)}] Any locally compact group is a hypergroup. See the monographs \cite{bh95,La23} for
more examples including double coset spaces, polynomial hypergroups and orbit hypergroups etc.,
as well as their applications.
\item[{\rm(iii)}] In the sequel, we always assume that $K$ is a hypergroup with a left
\emph{Haar measure}, that is, a non-negative Radon measure $\mu$ on $K$ such that, for any $x\in K$,
$\dz_x\ast\mu=\mu$. It has not been proved that every hypergroup has a Haar measure.
However, any commutative
hypergroup, compact hypergroup, discrete hypergroup, and double coset hypergroup admits a Haar measure;
see \cite{bh95,La23} for more examples.
\end{enumerate}
\end{remark}

Denote by $\mz$ the space of all $\mu$-measurable complex-valued functions on $K$.
Recall also that the convolution of any two functions $f,g\in\mz$ is defined by setting,
for any $x\in K$,
\begin{align}\label{2e1}
(f\ast g)(x):=\int_K f(y)g(y^-\ast x)\,d\mu(y),
\end{align}
where
$$g(y^-\ast x):=\int_K g(t)\,d(\dz_{y^-}\ast\dz_x)(t).$$
In addition, for any $x,y\in K$ and $U,V\subset K$, let
$$\{x\}\ast\{y\}:=\supp(\dz_x\ast\dz_y)\quad {\rm and}\quad
U\ast V:=\bigcup_{x\in U}\{x\}\ast V:=\bigcup_{x\in U}\bigcup_{y\in V}\{x\}\ast\{y\}.$$

\subsection{Orlicz and Calder\'{o}n--Lozanovski\u{\i} spaces}\label{s2.2}

To formulate the main results, we also need some notions on Orlicz spaces from \cite{kr61,rr91}.
Recall that a function $\vz:\ \mathbb{R}\to[0,\fz)$ is called a \emph{Young function}
if it is convex, even, left continuous on $(0,\fz)$ and $\vz(0)=0$; we also assume that $\vz$ is
neither identically zero nor identically infinite on $(0,\fz)$. Moreover, for any Young function
$\vz$, let
\begin{align*}
a_\vz:=\sup\lf\{x\in\mathbb{R}:\ \vz(x)=0\r\}\ {\rm and}\
b_\vz:=\sup\lf\{x\in\mathbb{R}:\ \vz(x)<\fz\r\}.
\end{align*}
Then it is continuous
on $[0, b_\vz)$, nondecreasing on $[0,\fz)$ and strictly increasing on $[a_\vz, b_\vz]$;
see \cite{agms}.

Let $K$ be a hypergroup with a left Haar measure $\mu$ and $\vz$ a Young function. For any $f\in\mz$,
define
$$\varrho_{\vz}(f):=\int_K\vz(|f(x)|)\,d\mu(x)$$
and
\begin{equation*}
\|f\|_{\lv}:=\inf\lf\{\lz\in(0,\fz):\ \varrho_{\vz}(f/\lz)\le1\r\},
\end{equation*}
which is called the \emph{Luxemburg norm} of $f$. Furthermore, the \emph{Orlicz space} $\lv$ is defined
to be the set of all $f\in\mz$ satisfying that $\varrho_{\vz}(f)<\fz$,
equipped with the Luxemburg norm $\|\cdot\|_{\lv}$.
Observe that the Orlicz space $\lv$ is a Banach space.

The following concept of Calder\'{o}n--Lozanovski\u{\i} spaces, which are defined in the similar way
as Orlicz spaces, comes from \cite{agms,Ca64,Lo65}.

\begin{definition}\label{2d2}
Let $K$ be a hypergroup with a left Haar measure $\mu$.
\begin{enumerate}
\item[(i)] A Banach space $\Oz=(\Oz,\|\cdot\|_\Oz)$ is called a \emph{Banach ideal space} on $K$
if $\Oz$ is a linear subspace of $\mz$ and
satisfies the following condition: if $f\in\Oz$, $g\in\mz$ and $|g(x)|\le|f(x)|$ for almost every
$x\in K$, then $g\in\Oz$ and $\|g\|_\Oz\le\|f\|_\Oz$.
\item[(ii)] For a given Banach ideal space $\Oz$ on $K$ and any Young function $\vz$,
the \emph{Calder\'{o}n--Lozanovski\u{\i} space} $\cl$ is defined
to be the set of all $f\in\mz$ satisfying that $I_{\vz}(f)<\fz$,
equipped with the \emph{Luxemburg norm}
\begin{equation*}
\|f\|_{\cl}:=\inf\lf\{\lz\in(0,\fz):\ I_{\vz}(f/\lz)\le1\r\},
\end{equation*}
where
\begin{align*}
I_{\vz}(f):=
\begin{cases}
\|\vz(|f|)\|_{\Oz}
&\text{if}\ \vz(|f|)\in\Oz,\\
\infty
&\text{otherwise}.
\end{cases}
\end{align*}
\end{enumerate}
\end{definition}

For any $p\in(0,\fz]$, the Lebesgue space $L^p(K)$ is defined to be the set of all $f\in\mz$
such that, when $p\in(0,\fz)$,
$$\|f\|_{L^p(K)}:=\lf[\int_K|f(x)|^p\,d\mu(x)\r]^{1/p}<\fz$$
and
$$\|f\|_{L^\fz(K)}:=\esup_{x\in K}|f(x)|<\fz.$$

\begin{remark}\label{2r2}
\begin{enumerate}
\item[{\rm(i)}] When $\Oz=L^1(K)$, then $\cl$ is the Orlicz space $\lv$.
When $\Oz$ is the Lorentz space on $K$, then $\cl$ becomes the corresponding Orlicz--Lorentz space.
\item[{\rm(ii)}] If $\vz(t)=t^p$ for any given $p\in[1,\fz)$, then $\cl$ is just the $p$
convexification $\Oz^{(p)}$ of $\Oz$ with the norm
$$\|\cdot\|_{\Oz^{(p)}}=\||\cdot|^p\|_{\Oz}^{1/p}.$$
In addition, if $\vz(t)=0$ for $t\in[0,1]$ and $\vz(t)=\infty$ otherwise, then $\cl$ goes back to
the space $L^\fz(K)$.
\end{enumerate}
\end{remark}

\subsection{Some assumptions on the Banach ideal $\Oz$}\label{s2.3}

Let $K$ be a hypergroup with a left Haar measure $\mu$.
Nonnegative (real) functions $f,\,g\in\mz$ are called \emph{equimeasurable}, if for any $s\in[0,\fz)$,
$$\mu(\{x\in K:\ |f(x)|>s\})=\mu(\{x\in K:\ |g(x)|>s\}).$$
In what follows, we always assume that $\Oz$ is a Banach ideal in $\mz$ which satisfies additionally
the conditions:
\begin{enumerate}
\item[{\rm(I)}] For any sequence of nonnegative (real) functions $\{f_n\}_{n\in\nn}\subset\Oz$
and $f\in\mz$ satisfying $f_n\uparrow f$ almost everywhere on $K$, if $f\in\Oz$, then
$\|f_n\|_\Oz\uparrow \|f\|_\Oz$ and, if $f\notin\Oz$, then $\|f_n\|_\Oz\to\fz$;
\item[{\rm(II)}] for any measurable subset $E\subset K$, if $\mu(E)<\fz$, then $\mathbf{1}_E\in\Oz$;
\item[{\rm(III)}] for any measurable subset $E\subset K$, if $\mu(E)<\fz$, then there exists a positive
constant $C_E$, depending on $E$, such that, for any $f\in \Oz$,
$$\int_E|f(x)|d\mu(x)\le C_E\|f\|_\Oz;$$
\item[{\rm(IV)}] if $f,\,g\in\Oz$ are equimeasurable, then $\|f\|_\Oz=\|g\|_\Oz$.
\end{enumerate}

\begin{remark}\label{2r3}
If $\Oz$ consists of real functions, then the above conditions mean that $\Oz$ is a Banach function space
(see \cite[p.\,3, Definition 1.3]{bs88}) and also a rearrangement-invariant space
(see \cite[p.\,59, Definition 4.1]{bs88}); if $\Oz$ is complex, then in such case, the real space
$(\Oz^{\mathbb{R}},\|\cdot\|_\Oz)$ is a Banach function space and also a rearrangement-invariant space,
where $\Oz^{\mathbb{R}}:=\{{\rm Re}(f):f\in\Oz\}$.
\end{remark}

Note that, for any measurable subsets $E,\,F\subset K$ with $\mu(E)=\mu(F)$, their characteristic functions
$\mathbf{1}_E$ and $\mathbf{1}_F$ are equimeasurable and hence $\|\mathbf{1}_E\|_\Oz=\|\mathbf{1}_F\|_\Oz$.
Therefore, we can define a function $\Phi_\Oz:\ [0,\fz)\to[0,\fz)$ as follows: for any measurable subset
$E\subset K$ with $\mu(E)<\fz$,
\begin{align}\label{2e2}
\Phi_\Oz(\mu(E)):=\|\mathbf{1}_E\|_\Oz,
\end{align}
which is called the \emph{fundamental function} of $\Oz$; see \cite[p.\,65, Definition 5.1]{bs88}.

Then \cite[p.\,67, Corollary 5.3]{bs88} implies the following Lemma \ref{2l1}.

\begin{lemma}\label{2l1}
The fundamental function $\Phi_\Oz$ defined in \eqref{2e2} has the following properties:
\begin{enumerate}
\item[{\rm(i)}] $\Phi_\Oz(0)=0$.
\item[{\rm(ii)}] $\Phi_\Oz$ is nondecreasing and continuous except (perhaps) at origin.
\item[{\rm(iii)}] The mapping $s\to\frac{\Phi_\Oz(s)}s$ is nonincreasing and, for any $\varepsilon\in(0,\fz)$,
there exists some $s_0\in(0,\fz)$ such that, for any $s\in[s_0,\fz)$,
$$\Phi_\Oz(s)\le\lf[\lim_{s\to\fz}\frac{\Phi_\Oz(s)}s+\varepsilon\r]s.$$
\end{enumerate}
\end{lemma}

Finally, we make two more assumptions on the Banach ideal $\Oz$:
\begin{enumerate}
\item[{\rm(V)}]
$\displaystyle\lim_{s\to0}\Phi_\Oz(s)=0$,
that is, the fundamental function $\Phi_\Oz$ is continuous at origin;
\item[{\rm(VI)}]
$\displaystyle\lim_{s\to\fz}\Phi_\Oz(s)=\fz$,
that is, the fundamental function $\Phi_\Oz$ is unbounded.
\end{enumerate}

\begin{remark}\label{2r4}
\begin{enumerate}
\item[{\rm(i)}] By \cite[p.\,67, Theorem 5.5]{bs88}, we find that, if $\Oz$ is separable,
then $\Phi_\Oz$ is continuous at origin.
\item[{\rm(ii)}] The assumption (VI) is also natural, for instance, it (and also the
assumption (V)) is satisfied when $K$ is not compact and $\Oz=L^p(K)$ with $\|\cdot\|_\Oz=\|\cdot\|_{L^p(K)}$.
\item[{\rm(iii)}] We point out that, under the above assumptions (I)-(VI), the Banach ideal $\Oz$
can be taken as many function spaces such as Lebesgue spaces, Orlicz spaces, etc.
\end{enumerate}
\end{remark}

\subsection{Main results}\label{s2.4}

We begin with some concepts which are later used to state our first result.

\begin{definition}\label{2d3}
Let $\vz_1$ and $\vz_2$ be two Young functions. The \emph{product Calder\'{o}n--Lozanovski\u{\i} space}
is defined as
$$\clz\times\clzz :=\lf\{(f,g):\ f\in\clz,\ g\in\clzz\r\}$$
with the norm
$$\|(f,g)\|_{\clz\times\clzz}:=\max\lf\{\|f\|_{\clz},\ \|g\|_{\clzz}\r\}<\fz,$$
where $\{\Oz_{\vz_i}(K)\}_{i=1,2}$ are Calder\'{o}n--Lozanovski\u{\i} spaces as in Definition \ref{2d2}(ii).
\end{definition}

A hypergroup $K$ is said to be satisfying a condition $(\divideontimes)$ if, for every compact neighborhood
$U$ of the identity element of $K$, there exist some $\omega\in(1,\fz)$ and
a strictly increasing sequence $\{n_k\}_{k\in\nn}\subset\nn$ such that, for any $k\in\nn$,
\begin{align}\label{2e3}
\mu(\overbrace{U\ast\cdots\ast U}^{2n_k\ \mathrm{times}})
<\omega\mu(\overbrace{U\ast\cdots\ast U}^{n_k\ \mathrm{times}}),
\end{align}
and there exists some $k_0\in\nn$ large enough such that, for any
$x\in(\overbrace{U\ast\cdots\ast U}^{n_{k_0}\ \mathrm{times}})$,
\begin{align}\label{2e3'}
\dz_{x}\ast\dz_{x^-}=\dz_e=\dz_{x^-}\ast\dz_{x}.
\end{align}

\begin{remark}\label{2r5}
By an argument similar to the proof of \cite[Proposition 16.28]{Pi84}, we know that every hypergroup
having the so-called polynomial growth (see \cite[Definition 2.5.11(b)]{bh95} for the precise definition)
satisfies the condition \eqref{2e3}; see also \cite[Chapters 2 and 3]{bh95} for examples of hypergroups
having polynomial growth.
\end{remark}

The \emph{inverse} of a Young function $\vz$ is defined by setting, for any $y\in[0,\fz)$,
$$\vz^{-1}(y):=\sup\{x\in[0,\fz):\ \vz(x)\le y\}.$$
Then the first main result of this article is stated as follows.

\begin{theorem}\label{t1}
Assume that $K$ is a non-compact hypergroup satisfying the condition $(\divideontimes)$.
Let $\{\vz_i\}_{i=1,2,3}$ be Young functions with $\vz_i(b_{\vz_i})>0$ for $i=1,2,3$ and
\begin{align}\label{2e4}
\liminf_{x\to0}\frac{\vz_1^{-1}(x)\vz_2^{-1}(x)}{x\vz_3^{-1}(x)}=\fz.
\end{align}
Then the set
$$E=\lf\{(f,g)\in\clz\times\clzz:\ |f|\ast|g|\in\clzzz\r\}$$
is of first category in $\clz\times\clzz$.
\end{theorem}

A hypergroup $K$ is said to be satisfying a condition $(\divideontimes\divideontimes)$ if,
for any compact neighborhood $U$ of the identity element of $K$, there exist some $\tau\in(1,\fz)$ and
a sequence of subsets $\{U_n\}_{n\in\nn}\subset U$ such that
$$\lim_{n\to\fz}\mu(U_n)=0,\quad \mu\lf(U_n^-\ast U_n\r)\le\tau\mu(U_n),\
\forall\,n\in\nn,$$
and there exists some $n_0\in\nn$ large enough such that, for any
$x\in U_{n_0}$,
\begin{align}\label{2e4'}
\dz_{x}\ast\dz_{x^-}=\dz_e=\dz_{x^-}\ast\dz_{x}.
\end{align}
Using this, we give the second main result of this article.

\begin{theorem}\label{t2}
Assume that $K$ is a compact hypergroup fulfilling the condition $(\divideontimes\divideontimes)$.
Let $\{\vz_i\}_{i=1,2,3}$ be Young functions with
\begin{align*}
\liminf_{x\to\fz}\frac{\vz_1^{-1}(x)\vz_2^{-1}(x)}{x\vz_3^{-1}(x)}=\fz.
\end{align*}
Then the set
$$\lf\{(f,g)\in\clz\times\clzz:\ |f|\ast|g|\in\clzzz\r\}$$
is of first category in $\clz\times\clzz$.
\end{theorem}

A hypergroup $K$ is said to be satisfying an adapted Leptin condition if, for every compact subset
$U\subset K$ and any $\varepsilon\in(0,\fz)$, there exists some compact subset $V\subset K$
such that
$$0<\mu(V)<\fz,\quad \mu(U\ast V)<(1+\varepsilon)\mu(V)$$
and, for any $x\in V$,
\begin{align}\label{2e5'}
\dz_{x}\ast\dz_{x^-}=\dz_e=\dz_{x^-}\ast\dz_{x}.
\end{align}
The following Theorem \ref{t3} is the third main result of this article.

\begin{theorem}\label{t3}
Assume that $K$ is a non-compact hypergroup satisfying the adapted Leptin condition.
Let $\vz,\,\phi$ be two Young functions such that $\vz(b_{\vz})>0$, $\phi(b_{\phi})=\fz$ and
\begin{align}\label{2e5}
\lim_{x\to0}\frac{\vz^{-1}(x)}{x}=\fz.
\end{align}
Then the set
$$E=\lf\{(f,g)\in\cl\times\clp:\ |f|\ast|g|\in\clp\r\}$$
is of first category in $\cl\times\clp$.
\end{theorem}

\begin{remark}\label{2r6}
As a special case of hypergroups, if $K$ is a locally compact group, then the conditions
\eqref{2e3'}, \eqref{2e4'} and \eqref{2e5'} are automatically satisfied. In this case,
\begin{enumerate}
\item[{\rm(i)}]
Theorem \ref{t1} comes back to \cite[Theorem 3.7]{agms} which is an extension of the first part
of \cite[Theorem 11]{km89}; see \cite[Remark 3.8]{agms} for more details.
\item[{\rm(ii)}]
Theorem \ref{t2} with $\Oz=L^1(K)$ is just \cite[Theorem 3.9]{agms} which is a topological strengthening
of \cite[Theorem 14]{km89} with a slightly weaker condition.
\item[{\rm(iii)}] The adapted Leptin condition used in Theorem \ref{t3} is coincident with the so-called
Leptin condition, which is equivalent to the amenability of a locally compact group. In this sense,
Theorem \ref{t3} extends \cite[Theorem 3.11]{agms} from locally compact groups to hypergroups.
We point out that the Leptin condition in Theorem \ref{t3} cannot be dropped even for locally compact groups;
see \cite{ks60} or \cite[Remark 3.12]{agms} for more details.
\end{enumerate}
\end{remark}

Recall that the \emph{modular function} on $K$ with the
Haar measure $\mu$ is a positive real valued function $\Delta$ on $K$ satisfying that, for any $x\in K$,
$$\mu\ast\dz_{x^-}=\Delta(x)\mu.$$
The hypergroup $K$ is called \emph{unimodular} if $\Delta\equiv1$;
see \cite[p.\,27]{Je75} or \cite[Definition 1.3.23]{bh95} for more details.

A compact set $U\subset K$ is called satisfying a condition $(\star)$
if there exists a compact symmetric neighborhood $V$ of the identity element $e$ in $K$ such that
$U\subset V$ and, for any $x\in V$ and any sequence $\{k_n\}_{n\in\nn}\subset K$ with $\Delta(k_n)\le1$,
it holds true that
\begin{align}\label{2e6}
\dz_{x}\ast\dz_{x^-}=\dz_e=\dz_{x^-}\ast\dz_{x}
\end{align}
and
\begin{align}\label{2e7}
\mu(\{k_n\}\ast V\ast V)\le C_{(U,V)}\mu(\{k_n\}\ast V),
\end{align}
where $C_{(U,V)}\in[1,\fz)$ is a constant depending on $U$ and $V$, but independent of $\{k_n\}_{n\in\nn}$.

Now, we state the fourth main result of this article.

\begin{theorem}\label{t4}
Let $K$ be a non-compact hypergroup.
Assume that $\vz,\,\phi$ are two Young functions such that $\vz(b_{\vz})>0$, $\phi(b_{\phi})>0$ and
\begin{align*}
\liminf_{x\to0}\frac{\vz^{-1}(x)\phi^{-1}(x)}{x}=\fz.
\end{align*}
Then, for each compact set $U$ satisfying $\mu(U)>0$ and the condition $(\star)$, the set
$$E_U=\lf\{(f,g)\in\cl\times\clp:\
(|f|\ast |g|)(x)\ is\ well\ defined\ at\ some\ point\ x\in U\r\}$$
is of first category in $\cl\times\clp$.
\end{theorem}

A compact set $U\subset K$ is called satisfying a condition $(\star\star)$
if there exist some $a\in K$ and a compact symmetric neighborhood $V$ of the identity element $e$ in $K$
such that $\Delta(a)<1$, $U\subset V$ and, for any $x\in V$ and $k\in\nn$,
it holds true that
\begin{align}\label{2e9}
\dz_{x}\ast\dz_{x^-}=\dz_e=\dz_{x^-}\ast\dz_{x},
\end{align}
\begin{align}\label{2e10}
\mu\lf(\{a\}^k\ast V\ast V\r)\le C_{(U,V,a)}\mu(V\ast V)
\end{align}
and
\begin{align}\label{2e11}
\lim_{k\to\fz}\mu\lf(V\ast \{a^-\}^k\r)=\fz,
\end{align}
where $C_{(U,V,a)}\in[1,\fz)$ is a constant depending on $U,V,a$, but independent of $k$,
and
$$\{a\}^k:=\overbrace{\{a\}\ast\cdots\ast \{a\}}^{k\ \mathrm{times}}.$$

Finally, we present the last main result of this article.

\begin{theorem}\label{t5}
Let $K$ be a non-compact hypergroup.
Assume that $\vz,\,\phi$ are two Young functions such that $\vz(b_{\vz})>0$, $\phi(b_{\phi})>0$ and
\begin{align}\label{2e12}
\liminf_{x\to0}\frac{\vz^{-1}(x)}{x}=\fz.
\end{align}
Then, for each compact set $U$ satisfying $\mu(U)>0$ and the condition $(\star\star)$, the set
$$E_U=\lf\{(f,g)\in\cl\times\clp:\
(|f|\ast |g|)(x)\ is\ well\ defined\ at\ some\ point\ x\in U\r\}$$
is of first category in $\cl\times\clp$.
\end{theorem}

\begin{remark}\label{2r7}
As a special cases of hypergroups:
\begin{enumerate}
\item[{\rm(i)}]
If $K$ is a locally compact group, then the condition $(\star)$ is satisfied with
$$C_{U,V}=\frac{\mu(V\ast V)}{\mu(V)}.$$
In this case, Theorem \ref{t4} comes back to \cite[Theorem 3.13]{agms}.
\item[{\rm(ii)}]
If $K$ is a non-unimodular locally compact group, then $K$ is non-compact and,
for every compact set $U\subset K$,
there exist some $a\in K$ and a compact symmetric neighborhood $V$ of the identity element $e$ in $K$
such that $\Delta(a)<1$, $U\subset V$ and, for any $x\in V$ and $k\in\nn$,
it holds true that $x\cdot x^-=e=x^-\cdot x$,
\begin{align*}
\mu\lf(a^kV^2\r)=\mu(V^2)
\end{align*}
and
\begin{align*}
\lim_{k\to\fz}\mu\lf(Va^{-k}\r)=\lim_{k\to\fz}\mu(V)\Delta\lf(a^{-k}\r)
=\lim_{k\to\fz}\mu(V)\lf[\Delta^{-}(a)\r]^k=\fz
\end{align*}
since $\mu(V)\ge\mu(U)>0$ and $\Delta^{-}(a)=\frac1{\Delta(a)}>1$. This illustrates that
the condition $(\star\star)$ is satisfied with $C_{U,V,a}=1.$
In this case, Theorem \ref{t5} is just \cite[Theorem 3.14]{agms}.
\end{enumerate}
\end{remark}

Compared with the corresponding conclusions for locally compact groups determined in \cite{agms},
the above results on hypergroups need some additional conditions. This also leads us to the following
natural questions, which are still unclear so far.

\setcounter{theorem}{0}
\renewcommand{\thetheorem}{}
\begin{question}
\begin{enumerate}
\item[{\rm(i)}]
Whether Theorems \ref{t1}, \ref{t2} and \ref{t3} hold true without, respectively,
the conditions \eqref{2e3'}, \eqref{2e4'} and \eqref{2e5'}.
\item[{\rm(ii)}]
Whether Theorem \ref{t4} without the condition $(\star)$ and Theorem \ref{t5}
without the condition $(\star\star)$ for non-unimodular hypergroups, respectively, hold true.
\end{enumerate}
\end{question}

\section{Proofs of Theorems \ref{t1} and \ref{t2}}\label{s3}

To show Theorems \ref{t1} and \ref{t2},
we need some technical lemmas. The following Lemmas \ref{3l1} and \ref{3l2}
come from \cite[Lemmas 2.1 and 2.3]{agms}, respectively.

\setcounter{theorem}{0}
\renewcommand{\thetheorem}{\arabic{section}.\arabic{theorem}}

\begin{lemma}\label{3l1}
For any Young function $\vz$ and its inverse $\vz^{-1}$, it holds true that
\begin{enumerate}
\item[{\rm(i)}] for any $x\in[0,\fz)$, $\vz(\vz^{-1}(x))\le x$;
\item[{\rm(ii)}] for any $x\in[0,\vz(b_\vz)]$, $\vz(\vz^{-1}(x))=x$;
\item[{\rm(iii)}] if $\vz(x)<\fz$, then $x\le\vz^{-1}(\vz(x))$;
\item[{\rm(iv)}] for any $x\in[a_\vz,b_\vz]$, $x=\vz^{-1}(\vz(x))$.
\end{enumerate}
\end{lemma}

\begin{lemma}\label{3l2}
Let $K$ be a hypergroup with a left Haar measure $\mu$ and $\vz$ a Young function.
Then, for any subset $E\subset K$ with $0<\mu(E)<\fz$,
\begin{align*}
\|\mathbf{1}_E\|_\cl=\lf[\vz^{-1}\lf(\frac1{\|\mathbf{1}_E\|_\Oz}\r)\r]^{-1}.
\end{align*}
\end{lemma}

For any $t\in[0,\fz)$, let
$$\Phi^{-1}_\Oz(t):=\sup\lf\{s\in[0,\fz):\ \Phi_\Oz(s)\le t\r\},$$
where $\Phi_\Oz$ denotes the fundamental function of $\Oz$ as in \eqref{2e2}.

Similarly to \cite[Lemma 3.6]{agms}, we easily obtain the succeeding Lemma \ref{3l3}.

\begin{lemma}\label{3l3}
\begin{enumerate}
\item[{\rm(i)}] The function $\Phi^{-1}_\Oz$ is nondecreasing and
    $\displaystyle\lim_{t\to0}\Phi^{-1}_\Oz(t)=0$;
\item[{\rm(ii)}] for any $s\in[0,\fz)$, $\Phi^{-1}_\Oz(\Phi_\Oz(s))\ge s$.
\end{enumerate}
\end{lemma}

Now, we give the proofs of Theorems \ref{t1} and \ref{t2}, respectively.

\begin{proof}[\textbf{Proof of Theorem \ref{t1}}]
For any $m\in\nn$, let
$$E_m:=\lf\{(f,g)\in\clz\times\clzz:\ \lf\|\,|f|\ast|g|\,\r\|_{\clzzz}< m\r\}.$$
Then $E=\bigcup_{m\in\nn}E_m$. To show Theorem \ref{t1}, it suffices to prove that, for each $m\in\nn$,
$E_m$ is nowhere dense.

To this end, for any given $D\in(0,\fz)$ and $m\in\nn$, by \eqref{2e4} and the continuity of every $\vz_i$
on $[0,\vz_i(b_{\vz_i}))$, we find that there exists
$$x_m\in\lf(0,\min\lf\{\vz_i\lf(b_{\vz_i}\r):\ i=1,2,3\r\}\r)$$
such that, for any $x\in(0,x_m]$,
\begin{align}\label{3e1}
\frac{D^2}{288(M+1)}\frac{\vz_1^{-1}(x)\vz_2^{-1}(x)}{x\vz_3^{-1}(x)}>m,
\end{align}
where $\displaystyle M:=\lim_{s\to\fz}\frac{\Phi_\Oz(s)}s$
and $\Phi_\Oz$ is the fundamental function as in \eqref{2e2}.
Observe that $K$ is non-compact. According to the assumption (VI) of $\Oz$ and Lemma \ref{2l1}(iii),
we can take a compact symmetric neighborhood $U$ of the identity element of $K$ satisfying that
\begin{align}\label{3e2}
\frac{1}{\|\mathbf{1}_U\|_{\Oz}}<x_m\quad {\rm and} \quad
\|\mathbf{1}_U\|_{\Oz}\le(M+1)\mu(U).
\end{align}
The condition $(\divideontimes)$ implies that, for such compact neighborhood $U$, there exist some $\omega\in(1,\fz)$ and
a strictly increasing sequence $\{n_k\}_{k\in\nn}\subset\nn$ such that
$\mu(U^{2n_k})<\omega\mu(U^{n_k})$ for any $k\in\nn$,
and there exists some $k_1\in\nn$ large enough such that, for any $x\in U^{n_{k_1}}$,
\begin{align}\label{3e2'}
\dz_{x}\ast\dz_{x^-}=\dz_e=\dz_{x^-}\ast\dz_{x},
\end{align}
here and thereafter, for any $k\in\nn$, let
$U^k:=\overbrace{U\ast\cdots\ast U}^{k\ \mathrm{times}}$.

Note that the sequence $\{\mu(U^{n_k})\}_{k\in\nn}$ is increasing by \cite[Lemma 3.3C]{Je75}. Thus, if
$\{\mu(U^{n_k})\}_{k\in\nn}$ is bounded, then it is convergent. From this, the fact that, for any $k\in\nn$,
$$\frac1{2\|\mathbf{1}_{U^{n_k}}\|_\Oz}\le\frac1{\|\mathbf{1}_{U}\|_\Oz}<x_m,$$
as well as \eqref{3e1} and \eqref{3e2}, it follows that there exists $k_2\in\nn$ such that,
for any $k\in[k_2,\fz)\cap\nn$, $\mu(U^{2n_k})<2\mu(U^{n_k})$ and
\begin{align}\label{3e3}
\frac{D^2\|\mathbf{1}_{U^{n_k}}\|_\Oz}{72(M+1)}\vz_1^{-1}\lf(\frac1{2\|\mathbf{1}_{U^{n_k}}\|_\Oz}\r)
\vz_2^{-1}\lf(\frac1{2\|\mathbf{1}_{U^{n_k}}\|_\Oz}\r)
>2m\vz_3^{-1}\lf(\frac1{2\|\mathbf{1}_{U^{n_k}}\|_\Oz}\r).
\end{align}
If $\displaystyle\lim_{k\to\fz}\mu(U^{n_k})=\fz$,
then, using \eqref{2e4} with $x=\frac1{\omega\|\mathbf{1}_{U^{n_k}}\|_\Oz}$, we can also choose $k_2\in\nn$
such that, for any $k\in[k_2,\fz)\cap\nn$,
\begin{align}\label{3e4}
\frac{D^2\|\mathbf{1}_{U^{n_k}}\|_\Oz}{72(M+1)}\vz_1^{-1}\lf(\frac1{\omega\|\mathbf{1}_{U^{n_k}}\|_\Oz}\r)
\vz_2^{-1}\lf(\frac1{\omega\|\mathbf{1}_{U^{n_k}}\|_\Oz}\r)
>\omega m\vz_3^{-1}\lf(\frac1{\omega\|\mathbf{1}_{U^{n_k}}\|_\Oz}\r).
\end{align}
Based on \eqref{3e3} and \eqref{3e4}, in the remainder of this proof, we proceed with the form of \eqref{3e4}
(all it has to do is set $\omega=2$ in the case when $\{\mu(U^{n_k})\}_{k\in\nn}$ is bounded).

Now, let $F_1:=U^{n_{k_0}}$ and $F_2:=U^{2n_{k_0}}$, where $k_0=\max\{k_1,k_2\}$. Then
$\mu(F_2)<\omega\mu(F_1)$ and, by the fact that $\omega>1$ and Lemma \ref{2l1}(iii), we have
$$\frac{\Phi_\Oz(\omega\mu(F_1))}{\omega\mu(F_1)}\le\frac{\Phi_\Oz(\mu(F_1))}{\mu(F_1)}.$$
This, combined with \eqref{2e2} and Lemma \ref{2l1}(ii), implies that
\begin{align}\label{3e5}
\lf\|\mathbf{1}_{F_2}\r\|_\Oz=\Phi_\Oz\lf(\mu(F_2)\r)\le\Phi_\Oz\lf(\omega\mu(F_1)\r)
\le\omega\Phi_\Oz\lf(\mu(F_1)\r)=\omega\lf\|\mathbf{1}_{F_1}\r\|_\Oz.
\end{align}
Observe that $\mu(F_1)>0$. From this and Lemma \ref{3l3}(i), we deduce that there exists a constant
$d\in(0,D/6)$ such that
\begin{align}\label{3e6}
L\Phi_\Oz^{-1}\lf(\frac{6d\|\mathbf{1}_{F_2}\|_\Oz}{D}\r)
+\Phi_\Oz^{-1}\lf(\frac{6d\|\mathbf{1}_{F_1}\|_\Oz}{D}\r)
\le\frac12\mu(F_1),
\end{align}
where $\displaystyle L:=\sup_{x\in F_2}\Delta(x^-)$ and $\Delta$ is the modular function with respect
to the Haar measure $\mu$ on $K$.

For any $(f,g)\in E_m$ and $x\in K$, let
\begin{align*}
\widetilde{f}(x):=
\begin{cases}
f(x)
&\text{if}\ x\notin F_1,\\
f(x)+C_f
&\text{if}\ x\in F_1\ \text{and}\ \text{Re}(f(x))\ge0,\\
f(x)-C_f
&\text{if}\ x\in F_1\ \text{and}\ \text{Re}(f(x))<0,
\end{cases}
\end{align*}
and
\begin{align*}
\widetilde{g}(x):=
\begin{cases}
g(x)
&\text{if}\ x\notin F_2,\\
g(x)+C_g
&\text{if}\ x\in F_2\ \text{and}\ \text{Re}(g(x))\ge0,\\
g(x)-C_g
&\text{if}\ x\in F_2\ \text{and}\ \text{Re}(g(x))<0,
\end{cases}
\end{align*}
where
\begin{align}\label{3e7}
C_f:=\frac{D}{3\|\mathbf{1}_{F_1}\|_{\clz}}
\xlongequal{\text{Lem.}\ \ref{3l2}}\frac{D}3\vz_1^{-1}\lf(\frac1{\|\mathbf{1}_{F_1}\|_\Oz}\r)
\end{align}
and
\begin{align}\label{3e8}
C_g:=\frac{D}{3\|\mathbf{1}_{F_2}\|_{\clzz}}
\xlongequal{\text{Lem.}\ \ref{3l2}}\frac{D}3\vz_2^{-1}\lf(\frac1{\|\mathbf{1}_{F_2}\|_\Oz}\r).
\end{align}
Then
$$\lf\|f-\widetilde{f}\,\r\|_{\clz}=\lf\|C_f\mathbf{1}_{F_1}\r\|_{\clz}=\frac D3<D$$
and
$$\lf\|g-\widetilde{g}\,\r\|_{\clzz}=\lf\|C_g\mathbf{1}_{F_2}\r\|_{\clzz}=\frac D3<D.$$
Using this, Definition \ref{2d3} and the fact that $d\in(0,D/6)$, we further conclude that
$$B\lf(\lf(\widetilde{f},\,\widetilde{g}\r),d\r)\subset B\lf((f,g),D\r).$$
Therefore, to prove that $E_m$ is nowhere dense for each $m\in\nn$,
it suffices to show that $B((\widetilde{f},\,\widetilde{g}),d)\cap E_m=\emptyset$.

To achieve this, for any $(u,v)\in B((\widetilde{f},\,\widetilde{g}),d)$, let
\begin{align}\label{3e7'}
F_{1,u}:=\lf\{x\in F_1:\ |u(x)|\le\frac{C_f}2\r\}
\end{align}
and
\begin{align}\label{3e8'}
F_{2,v}:=\lf\{x\in F_2:\ |v(x)|\le\frac{C_g}2\r\}.
\end{align}
Note that, according to the definition of complex moduli, we have $|\widetilde{f}(x)|\ge C_f$ for any
$x\in F_1$. Thus, by Lemma \ref{3l2}, \eqref{3e7'} and Definition \ref{2d3}, it is easy to see that
\begin{align*}
\vz_1^{-1}\lf(\frac1{\|\mathbf{1}_{F_{1,u}}\|_\Oz}\r)&=\frac{1}{\|\mathbf{1}_{F_{1,u}}\|_{\clz}}\\
&\ge\frac{C_f}{2\|(u-\widetilde{f}\ )\mathbf{1}_{F_{1,u}}\|_{\clz}}\\
&\ge\frac{C_f}{2\|(u-\widetilde{f}\ )\|_{\clz}}\ge\frac{C_f}{2d},
\end{align*}
which, together with Lemma \ref{3l1}(i), implies that
\begin{align}\label{3e9}
\vz_1\lf(\frac{C_f}{2d}\r)\le\vz_1\lf(\vz_1^{-1}\lf(\frac1{\|\mathbf{1}_{F_{1,u}}\|_\Oz}\r)\r)
\le\frac1{\|\mathbf{1}_{F_{1,u}}\|_\Oz}.
\end{align}
In addition, by Lemma \ref{3l1}(ii), the convexity of $\vz_1$ and \eqref{3e7}, we find that
\begin{align*}
\frac1{\|\mathbf{1}_{F_{1}}\|_\Oz}
&=\vz_1\lf(\vz_1^{-1}\lf(\frac1{\|\mathbf{1}_{F_{1}}\|_\Oz}\r)\r)\\
&\le\frac{6d}D\vz_1\lf(\frac D{6d}\vz_1^{-1}\lf(\frac1{\|\mathbf{1}_{F_{1}}\|_\Oz}\r)\r)\\
&=\frac{6d}D\vz_1\lf(\frac{C_f}{2d}\r).
\end{align*}
From this and \eqref{3e9}, it follows that
\begin{align}\label{3e10}
\frac1{\|\mathbf{1}_{F_{1}}\|_\Oz}
\le\frac{6d}D\frac1{\|\mathbf{1}_{F_{1,u}}\|_\Oz},
\end{align}
and similarly
\begin{align}\label{3e11}
\frac1{\|\mathbf{1}_{F_{2}}\|_\Oz}
\le\frac{6d}D\frac1{\|\mathbf{1}_{F_{2,v}}\|_\Oz}.
\end{align}

Observe that, by \cite[Theorem 5.3B]{Je75}, we know that
\begin{align*}
\mu\lf(\lf(F_{2,v}\r)^{-}\r)&=\mu^{-}\lf(F_{2,v}\r)=\int_{F_{2,v}}\,d\mu^{-}(x)\\
&=\int_{F_{2,v}}\Delta^{-}(x)\,d\mu(x)=\int_{F_{2,v}}\Delta(x^{-})\,d\mu(x)\\
&\le L\int_{F_{2,v}}\,d\mu(x)=L\mu\lf(F_{2,v}\r).
\end{align*}
By this, Lemma \ref{3l3}(ii), \eqref{2e2}, \eqref{3e10}, \eqref{3e11} and Lemma \ref{3l3}(i), we obtain
\begin{align*}
\mu\lf(\lf(F_{2,v}\r)^{-}\r)+\mu\lf(F_{1,u}\r)
&\le L\mu\lf(F_{2,v}\r)+\mu\lf(F_{1,u}\r)\\
&\le L\Phi^{-1}_\Oz\lf(\Phi_\Oz\lf(\mu\lf(F_{2,v}\r)\r)\r)
+\Phi^{-1}_\Oz\lf(\Phi_\Oz\lf(\mu\lf(F_{1,u}\r)\r)\r)\\
&=L\Phi^{-1}_\Oz\lf(\lf\|\mathbf{1}_{F_{2,v}}\r\|_\Oz\r)
+\Phi^{-1}_\Oz\lf(\lf\|\mathbf{1}_{F_{1,u}}\r\|_\Oz\r)\\
&\le L\Phi_\Oz^{-1}\lf(\frac{6d\|\mathbf{1}_{F_2}\|_\Oz}{D}\r)
+\Phi_\Oz^{-1}\lf(\frac{6d\|\mathbf{1}_{F_1}\|_\Oz}{D}\r),
\end{align*}
which, combined with \eqref{3e6}, gives us
\begin{align}\label{3e12}
\mu(F_1)-\mu\lf(\lf(F_{2,v}\r)^{-}\r)-\mu\lf(F_{1,u}\r)
\ge\mu(F_1)-\frac12\mu(F_1)=\frac12\mu(F_1).
\end{align}

For any $y\in F_1$, let
$$Y_y:=\lf(F_1\setminus F_{1,u}\r)\cap\lf(\{y\}\ast\lf(F_2\setminus F_{2,v}\r)^-\r).$$
Then
$$Y_y^-\ast\{y\}\subset F_2\setminus F_{2,v}.$$
Moreover, by \cite[Lemma 3.3C]{Je75}, \eqref{3e2'} and \cite[(2.1)]{bkt22},
we have
\begin{align}\label{3e13}
&\mu(F_1)-\mu\lf(F_{1,u}\r)-\mu\lf(\lf(F_{2,v}\r)^{-}\r)\\
&\quad=\mu\lf(F_1\setminus F_{1,u}\r)-\mu\lf(F_2^-\setminus\lf(F_2\setminus F_{2,v}\r)^{-}\r)\noz\\
&\quad=\mu\lf(\{y^-\}\ast\lf(F_1\setminus F_{1,u}\r)\r)
-\mu\lf(\lf(F_1\ast F_1\r)\setminus\lf(F_2\setminus F_{2,v}\r)^{-}\r)\noz\\
&\quad\le\mu\lf(\{y^-\}\ast\lf(F_1\setminus F_{1,u}\r)\r)
-\mu\lf(\lf(\{y^-\}\ast\lf(F_1\setminus F_{1,u}\r)\r)\setminus\lf(F_2\setminus F_{2,v}\r)^{-}\r)\noz\\
&\quad=\mu\lf(\lf(\{y^-\}\ast\lf(F_1\setminus F_{1,u}\r)\r)\cap\lf(F_2\setminus F_{2,v}\r)^{-}\r)\noz\\
&\quad=\mu\lf(\{y\}\ast\lf(\lf(\{y^-\}\ast\lf(F_1\setminus F_{1,u}\r)\r)
\cap\lf(F_2\setminus F_{2,v}\r)^{-}\r)\r)\noz\\
&\quad=\mu\lf(\lf(F_1\setminus F_{1,u}\r)
\cap\lf(\{y\}\ast\lf(F_2\setminus F_{2,v}\r)^{-}\r)\r)\noz\\
&\quad=\mu\lf(Y_y\r).\noz
\end{align}
Note that $Y_y\subset F_1\setminus F_{1,u}$ and, for any $x\in Y_y$,
$$\{x^-\}\ast\{y\}\subset Y_y^-\ast\{y\}\subset F_2\setminus F_{2,v}.$$
By this, \eqref{2e1}, Definition \ref{2d1}(iii), \eqref{3e7'} and \eqref{3e8'}, we conclude that,
for any $(u,v)\in B((\widetilde{f},\,\widetilde{g}),d)$ and $y\in F_1$,
\begin{align*}
\lf(|u|\ast|v|\r)(y)&=\int_K|u|(x)|v|(x^-\ast y)\,d\mu(x)\\
&=\int_{K}|u|(x)\int_K|v|(z)\,d(\dz_{x^-}\ast\dz_y)(z)\,d\mu(x)\\
&\ge\int_{Y_y}|u|(x)\int_{\{x^-\}\ast \{y\}}|v|(z)\,d(\dz_{x^-}\ast\dz_y)(z)\,d\mu(x)\\
&>\frac{C_fC_g}4\mu\lf(Y_y\r).
\end{align*}
From this, \eqref{3e7}, \eqref{3e8}, \eqref{3e12}, \eqref{3e13} and \eqref{3e5}, we deduce that
\begin{align}\label{3e14}
\lf(|u|\ast|v|\r)(y)
&>\frac{D^2\mu(F_1)}{72}\vz_1^{-1}\lf(\frac1{\|\mathbf{1}_{F_{1}}\|_\Oz}\r)
\vz_2^{-1}\lf(\frac1{\|\mathbf{1}_{F_{2}}\|_\Oz}\r)\\
&\ge\frac{D^2\mu(F_1)}{72}\vz_1^{-1}\lf(\frac1{\omega\|\mathbf{1}_{F_{1}}\|_\Oz}\r)
\vz_2^{-1}\lf(\frac1{\omega\|\mathbf{1}_{F_{1}}\|_\Oz}\r).\noz
\end{align}
In addition, by \eqref{2e2}, Lemma \ref{2l1}(iii) and \eqref{3e2}, we find that
$$\frac{\mu(F_1)}{\|\mathbf{1}_{F_{1}}\|_\Oz}=\frac{\mu(F_1)}{\Phi_\Oz(\mu(F_1))}
\ge\frac{\mu(U)}{\Phi_\Oz(\mu(U))}
=\frac{\mu(U)}{\|\mathbf{1}_{U}\|_\Oz}\ge\frac1{M+1}.$$
This, together with \eqref{3e14}, \eqref{3e3} and \eqref{3e4}, further implies that,
for any $(u,v)\in B((\widetilde{f},\,\widetilde{g}),d)$ and $y\in F_1$,
\begin{align*}
\lf(|u|\ast|v|\r)(y)
&>\frac{D^2\|\mathbf{1}_{F_{1}}\|_\Oz}{72(M+1)}\vz_1^{-1}\lf(\frac1{\omega\|\mathbf{1}_{F_{1}}\|_\Oz}\r)
\vz_2^{-1}\lf(\frac1{\omega\|\mathbf{1}_{F_{1}}\|_\Oz}\r)\frac{(M+1)\mu(F_1)}{\|\mathbf{1}_{F_{1}}\|_\Oz}\\
&\ge\omega m\vz_3^{-1}\lf(\frac1{\omega\|\mathbf{1}_{F_1}\|_\Oz}\r).
\end{align*}
Therefore, by the concavity of $\vz_3^{-1}$ and Lemma \ref{3l1}(ii), we know that,
for any $(u,v)\in B((\widetilde{f},\,\widetilde{g}),d)$,
\begin{align*}
\lf\|\vz_3\lf(\frac{|u|\ast|v|}m\r)\r\|_\Oz
&\ge\lf\|\vz_3\lf(\frac{|u|\ast|v|}m\r)\mathbf{1}_{F_1}\r\|_\Oz\\
&\ge\lf\|\vz_3\lf(\omega \vz_3^{-1}\lf(\frac1{\omega\|\mathbf{1}_{F_1}\|_\Oz}\r)\r)\mathbf{1}_{F_1}\r\|_\Oz\\
&\ge\lf\|\vz_3\lf(\vz_3^{-1}\lf(\frac1{\|\mathbf{1}_{F_1}\|_\Oz}\r)\r)\mathbf{1}_{F_1}\r\|_\Oz=1,
\end{align*}
that is, $\|\,|u|\ast|v|\,\|_{\clzzz}\ge m$ and hence $(u,v)\notin E_m$. This finishes the proof of
Theorem \ref{t1}.
\end{proof}

\begin{proof}[\textbf{Proof of Theorem \ref{t2}}]
Following the argument used in the proof of Theorem \ref{t1} with some slight modifications,
it is easy to see that Theorem \ref{t2} is valid; we omit the details. In particular, we should take
here $F_1:=U_{n_0}$ and $F_2:=U_{n_0}^-\ast U_{n_0}$ for some sufficiently large $n_0\in\nn$.
\end{proof}

\section{Proofs of Theorems \ref{t3}, \ref{t4} and \ref{t5}}\label{s4}

We first prove Theorem \ref{t3}.

\begin{proof}[\textbf{Proof of Theorem \ref{t3}}]
For any $m\in\nn$, define
$$E_m:=\lf\{(f,g)\in\cl\times\clp:\ \lf\|\,|f|\ast|g|\,\r\|_{\clp}< m\r\}.$$
Obviously, $E=\bigcup_{m\in\nn}E_m$. The proof will be completed if we show that, for every $m\in\nn$,
$E_m$ is nowhere dense.

For this purpose, for any given $D\in(0,\fz)$ and $m\in\nn$, by the non-compactness is of $K$,
the assumption (VI) of $\Oz$, Lemma \ref{2l1}(iii) and \eqref{2e5},
we can choose a large enough compact symmetric neighborhood $U$ of the identity element of $K$
such that
\begin{align}\label{4e1}
\frac{1}{\|\mathbf{1}_U\|_{\Oz}}<\vz(b_{\vz}),\quad
\|\mathbf{1}_U\|_{\Oz}=\Phi_\Oz(\mu(U))\le(M+1)\mu(U)
\end{align}
and
\begin{align}\label{4e2}
\frac{D^2\|\mathbf{1}_{U}\|_\Oz}{72(M+1)}\vz^{-1}\lf(\frac1{\|\mathbf{1}_{U}\|_\Oz}\r)
>2m,
\end{align}
where $\displaystyle M:=\lim_{s\to\fz}\frac{\Phi_\Oz(s)}s$ and $\Phi_\Oz$ is the fundamental function
as in \eqref{2e2}. Note that $K$ satisfies the adapted Leptin condition. Thus, there exists a compact
subset $V\in K$ such that
$$0<\mu(V)<\fz,\quad \mu(U\ast V)<2\mu(V)$$
and, for any $x\in V$,
\begin{align}\label{4e3}
\dz_{x}\ast\dz_{x^-}=\dz_e=\dz_{x^-}\ast\dz_{x}.
\end{align}
Let $F_2:=U\ast V$. Then, similarly to \eqref{3e5} and \eqref{3e6}, we have
\begin{align}\label{4e4}
\lf\|\mathbf{1}_{F_2}\r\|_\Oz=\Phi_\Oz\lf(\mu(F_2)\r)\le\Phi_\Oz\lf(2\mu(V)\r)
\le2\Phi_\Oz\lf(\mu(V)\r)=2\lf\|\mathbf{1}_{V}\r\|_\Oz
\end{align}
and
\begin{align*}
L\Phi_\Oz^{-1}\lf(\frac{6d\|\mathbf{1}_{F_2}\|_\Oz}{D}\r)
+\Phi_\Oz^{-1}\lf(\frac{6d\|\mathbf{1}_{U}\|_\Oz}{D}\r)
\le\frac12\mu(U),
\end{align*}
where $\displaystyle L:=\sup_{x\in F_2}\Delta(x^-)$ and $d\in(0,D/6)$.

Now, for any $(f,g)\in E_m$, define $\widetilde{f},\,\widetilde{g}$ on $K$ and the constants
$C_f,\,C_g$ as in the proof of Theorem \ref{t1} (here for the sets $U,\,F_2$ and functions
$\vz,\,\phi$, respectively). Then
$$\lf\|f-\widetilde{f}\,\r\|_{\cl}=\frac D3=\lf\|g-\widetilde{g}\,\r\|_{\clp}$$
and hence
$$B\lf(\lf(\widetilde{f},\,\widetilde{g}\r),d\r)\subset B\lf((f,g),D\r).$$
Therefore, we only need to verify that $B((\widetilde{f},\,\widetilde{g}),d)\cap E_m=\emptyset$.

To this end, for any $(u,v)\in B((\widetilde{f},\,\widetilde{g}),d)$, set
\begin{align*}
U_{u}:=\lf\{x\in U:\ |u(x)|\le\frac{C_f}2\r\}
\end{align*}
and
\begin{align*}
F_{2,v}:=\lf\{x\in F_2:\ |v(x)|\le\frac{C_g}2\r\}.
\end{align*}
Repeating the proof of \eqref{3e12} with some slight modifications, we obtain
\begin{align*}
\mu(U)-\mu\lf(\lf(F_{2,v}\r)^{-}\r)-\mu\lf(U_{u}\r)
\ge\frac12\mu(U).
\end{align*}

For any $y\in V$, let
$$Y_y:=\lf(U\setminus U_{u}\r)\cap\lf(\{y\}\ast\lf(F_2\setminus F_{2,v}\r)^-\r).$$
Then
$$Y_y\subset U\setminus U_{u}\quad {\rm and} \quad Y_y^-\ast\{y\}\subset F_2\setminus F_{2,v}.$$
In addition, from the symmetry of $U$ and the definition of $F_2$, it follows that
$$\{y^-\}\ast(U\setminus U_{u})\subset\{y^-\}\ast U=\{y^-\}\ast U^-\subset V^-\ast U^-=F_2^-.$$
This, combined with \eqref{4e3} and \cite[(2.1)]{bkt22}, gives us
\begin{align*}
&\mu(U)-\mu\lf(U_{u}\r)-\mu\lf(\lf(F_{2,v}\r)^{-}\r)\\
&\quad=\mu\lf(U\setminus U_{u}\r)-\mu\lf(F_2^-\setminus\lf(F_2\setminus F_{2,v}\r)^{-}\r)\\
&\quad\le\mu\lf(\{y^-\}\ast\lf(U\setminus U_{u}\r)\r)
-\mu\lf(\lf(\{y^-\}\ast\lf(U\setminus U_{u}\r)\r)\setminus\lf(F_2\setminus F_{2,v}\r)^{-}\r)\\
&\quad=\mu\lf(\lf(\{y^-\}\ast\lf(U\setminus U_{u}\r)\r)\cap\lf(F_2\setminus F_{2,v}\r)^{-}\r)\\
&\quad=\mu\lf(\{y\}\ast\lf(\lf(\{y^-\}\ast\lf(U\setminus U_{u}\r)\r)
\cap\lf(F_2\setminus F_{2,v}\r)^{-}\r)\r)\\
&\quad=\mu\lf(\lf(U\setminus U_{u}\r)
\cap\lf(\{y\}\ast\lf(F_2\setminus F_{2,v}\r)^{-}\r)\r)\\
&\quad=\mu\lf(Y_y\r).
\end{align*}
Via this and a calculation parallel to \eqref{3e14}, we conclude that
\begin{align*}
\lf(|u|\ast|v|\r)(y)
&>\frac{D^2\mu(U)}{72}\vz^{-1}\lf(\frac1{\|\mathbf{1}_{U}\|_\Oz}\r)
\phi^{-1}\lf(\frac1{\|\mathbf{1}_{F_{2}}\|_\Oz}\r).
\end{align*}
From this, \eqref{4e4}, \eqref{4e2} and \eqref{4e1}, we deduce that,
for any $(u,v)\in B((\widetilde{f},\,\widetilde{g}),d)$ and $y\in V$,
\begin{align*}
\lf(|u|\ast|v|\r)(y)
&>\frac{D^2\mu(U)}{72}\vz^{-1}\lf(\frac1{\|\mathbf{1}_{U}\|_\Oz}\r)
\phi^{-1}\lf(\frac1{2\|\mathbf{1}_{V}\|_\Oz}\r)\\
&=\frac{D^2\|\mathbf{1}_{U}\|_\Oz}{72(M+1)}\vz^{-1}\lf(\frac1{\|\mathbf{1}_{U}\|_\Oz}\r)
\phi^{-1}\lf(\frac1{2\|\mathbf{1}_{V}\|_\Oz}\r)\frac{(M+1)\mu(U)}{\|\mathbf{1}_{U}\|_\Oz}\\
&>2m\phi^{-1}\lf(\frac1{2\|\mathbf{1}_{V}\|_\Oz}\r).
\end{align*}
Thus, for any $(u,v)\in B((\widetilde{f},\,\widetilde{g}),d)$,
\begin{align*}
\lf\|\phi\lf(\frac{|u|\ast|v|}m\r)\r\|_\Oz
&\ge\lf\|\phi\lf(\frac{|u|\ast|v|}m\r)\mathbf{1}_{V}\r\|_\Oz\\
&\ge\lf\|\phi\lf(2\phi^{-1}\lf(\frac1{2\|\mathbf{1}_{V}\|_\Oz}\r)\r)\mathbf{1}_{V}\r\|_\Oz\\
&\ge\lf\|\phi\lf(\phi^{-1}\lf(\frac1{\|\mathbf{1}_{V}\|_\Oz}\r)\r)\mathbf{1}_{V}\r\|_\Oz=1
\end{align*}
and hence $(u,v)\notin E_m$, which completes the proof of Theorem \ref{t3}.
\end{proof}

Next, we show Theorem \ref{t4}.

\begin{proof}[\textbf{Proof of Theorem \ref{t4}}]
For any $m\in\nn$, let
$$E_{U,m}:=\lf\{(f,g)\in\cl\times\clp:\ \exists\,x\in U,\
\int_K |f|(y)|g|(y^-\ast x)\,d\mu(y)<m\r\}.$$
Then
$$E_U=\bigcup_{m\in\nn}E_{U,m}.$$
Therefore, to show Theorem \ref{t4}, it suffices to prove that, for each $m\in\nn$,
$E_{U,m}$ is nowhere dense.

Observe that the compact set $U$ satisfies the condition $(\star)$. Then there exists a compact
symmetric neighborhood $V$ of the identity element $e$ in $K$ such that $U\subset V$ and, \eqref{2e6}
and \eqref{2e7} hold true. Since $K$ is not compact, we can take an element $k\in K\setminus V$ such
that $k^-\in K\setminus V$ by the symmetry of $V$. By \cite[Theorem 5.3B]{Je75}, we have
$$\Delta(k)\Delta(k^-)=\Delta(k)\Delta^-(k)=1$$
and hence $0<\Delta(k)\le1$ or $0<\Delta(k^-)\le1$. This, together with the continuity of $\Delta$,
implies that there exists some $k_1\in K\setminus V$ such that $0<\Delta(k_1)\le1$. Let
$$V_1:=(\{k_1\}\ast V\ast V\ast V\ast V)\cup\lf( V\ast V\ast V\ast V\ast\{k_1^-\}\r).$$
Then, by \cite[Lemmas 3.2B and 3.4C]{Je75}, we find that $V_1$ is a compact symmetric subset of $K$.
Repeating the
above argument for $V_1$, we can also pick an element $k_2\in K\setminus V_1$ such that
$0<\Delta(k_2)\le1$. We claim that
\begin{align}\label{4e5}
\lf(V\ast\{k_1^-\}\r)\cap\lf(V\ast\{k_2^-\}\r)=\emptyset.
\end{align}
Indeed, if there exists an element $a\in(V\ast\{k_1^-\})\cap(V\ast\{k_2^-\})$, then
$e\in\{a^-\}\ast V\ast\{k_2^-\}$ by \cite[Lemma 4.1A]{Je75} and hence $e\in\supp(\dz_{b}\ast\dz_{k_2^-})$
with some $b\in \{a^-\}\ast V$. Thus, according to Definition \ref{2d1}(v) and \cite[Lemma 3.4C]{Je75},
we obtain
\begin{align*}
k_2=b\in\{a^-\}\ast V&\subset\{k_1\}\ast V^-\ast V\\
&=\{k_1\}\ast V\ast V\subset\{k_1\}\ast V\ast V\ast V\ast V.
\end{align*}
This contradicts the fact $k_2\in K\setminus V_1$, that is, \eqref{4e5} is true. Similarly,
we also have
\begin{align*}
\lf(\{k_1\}\ast V\ast V\r)\cap\lf(\{k_2\}\ast V\ast V\r)=\emptyset.
\end{align*}
Inductively, we get a sequence $\{k_n\}_{n\in\nn}\subset K$ such that $0<\Delta(k_n)\le1$,
$$k_n\in K\setminus \bigcup_{i\in[1,n-1]\cap\nn}V_{i},\quad\forall\,n\in[2,\fz)\cap\nn,$$
where
$$V_i:=\lf(\{k_i\}\ast V\ast V\ast V\ast V\r)\cup\lf(V\ast V\ast V\ast V\ast\{k_i^-\}\r)$$
and, for any $\ell,n\in\nn$ with $\ell\neq n$,
\begin{align}\label{4e6}
\lf(V\ast\{k_\ell^-\}\r)\cap\lf(V\ast\{k_n^-\}\r)=\emptyset
\end{align}
and
\begin{align}\label{4e7}
\lf(\{k_\ell\}\ast V\ast V\r)\cap\lf(\{k_n\}\ast V\ast V\r)=\emptyset.
\end{align}

For any given $D\in(0,\fz)$ and $m\in\nn$, from the assumption $\mu(V)\ge\mu(U)>0$ and Lemma \ref{2l1}(iii),
we infer that there exists some $N\in\nn$ large enough such that the following two sets
$$F_1=\bigcup_{n\in[1,N]\cap\nn}\lf(\{k_n\}\ast V\r)\quad{\rm and}\quad
F_2=\bigcup_{n\in[1,N]\cap\nn}\lf(\{k_n\}\ast V\ast V\r)$$
satisfy
\begin{align*}
\frac{1}{\|\mathbf{1}_{F_1}\|_{\Oz}}<\vz(b_{\vz}),\quad
\frac{1}{\|\mathbf{1}_{F_2}\|_{\Oz}}<\phi(b_{\phi}),
\end{align*}
\begin{align}\label{4e9}
\|\mathbf{1}_{F_1}\|_{\Oz}=\Phi_\Oz(\mu(F_1))\le(M+1)\mu(F_1)
\end{align}
and
\begin{align}\label{4e10}
\frac{D^2\|\mathbf{1}_{F_1}\|_\Oz}{72(M+1)}\vz^{-1}\lf(\frac1{C_{(U,V)}\|\mathbf{1}_{F_1}\|_\Oz}\r)
\phi^{-1}\lf(\frac1{C_{(U,V)}\|\mathbf{1}_{F_1}\|_\Oz}\r)>m,
\end{align}
where $\displaystyle M:=\lim_{s\to\fz}\frac{\Phi_\Oz(s)}s$ and $C_{(U,V)}$ is as in \eqref{2e7}.
In addition, by \eqref{4e7}, \eqref{2e7} and \eqref{4e6}, we see that
\begin{align*}
\mu(F_2)&=\sum_{n\in[1,N]\cap\nn}\mu\lf(\{k_n\}\ast V\ast V\r)\\
&\le C_{(U,V)}\sum_{n\in[1,N]\cap\nn}\mu(\{k_n\}\ast V)=C_{(U,V)}\mu(F_1).
\end{align*}
This, combined with \eqref{2e2} and (iii) and (ii) of Lemma \ref{2l1}, implies that
\begin{align*}
\lf\|\mathbf{1}_{F_2}\r\|_\Oz&=\Phi_\Oz\lf(\mu(F_2)\r)\le\Phi_\Oz\lf(C_{(U,V)}\mu(F_1)\r)\\
&\le C_{(U,V)}\Phi_\Oz\lf(\mu(F_1)\r)=C_{(U,V)}\lf\|\mathbf{1}_{F_1}\r\|_\Oz.
\end{align*}
Using this, \eqref{4e9} and \eqref{4e10}, we conclude that
\begin{align}\label{4e11}
\frac{D^2\mu(F_1)}{72}\vz^{-1}\lf(\frac1{\|\mathbf{1}_{F_1}\|_\Oz}\r)
\phi^{-1}\lf(\frac1{\|\mathbf{1}_{F_2}\|_\Oz}\r)>m.
\end{align}
Moreover, similarly to \eqref{3e6}, we have
\begin{align*}
L\Phi_\Oz^{-1}\lf(\frac{6d\|\mathbf{1}_{F_2}\|_\Oz}{D}\r)
+\Phi_\Oz^{-1}\lf(\frac{6d\|\mathbf{1}_{F_1}\|_\Oz}{D}\r)
\le\frac12\mu(F_1),
\end{align*}
where $\displaystyle L:=\sup_{x\in F_2}\Delta(x^-)$ and $d\in(0,D/6)$.

Now, for any $(f,g)\in E_{U,m}$, define $\widetilde{f},\,\widetilde{g}$ on $K$ and the constants
$C_f,\,C_g$ as in the proof of Theorem \ref{t1} (here for the sets $F_1,\,F_2$ and functions
$\vz,\,\phi$, respectively). Then
$$\lf\|f-\widetilde{f}\,\r\|_{\cl}=\frac D3=\lf\|g-\widetilde{g}\,\r\|_{\clp}$$
and hence
$$B\lf(\lf(\widetilde{f},\,\widetilde{g}\r),d\r)\subset B\lf((f,g),D\r).$$
Clearly, it suffices to show that $B((\widetilde{f},\,\widetilde{g}),d)\cap E_{U,m}=\emptyset$.

For this purpose, for any $(u,v)\in B((\widetilde{f},\,\widetilde{g}),d)$, let
\begin{align*}
F_{1,u}:=\lf\{x\in F_1:\ |u(x)|\le\frac{C_f}2\r\}
\end{align*}
and
\begin{align*}
F_{2,v}:=\lf\{x\in F_2:\ |v(x)|\le\frac{C_g}2\r\}.
\end{align*}
Then, by a calculation parallel to \eqref{3e12}, we conclude that
\begin{align*}
\mu(F_1)-\mu\lf(\lf(F_{2,v}\r)^{-}\r)-\mu\lf(F_{1,u}\r)
\ge\frac12\mu(F_1).
\end{align*}
For any $y\in V$, let
$$Y_y:=\lf(F_1\setminus F_{1,u}\r)\cap\lf(\{y\}\ast\lf(F_2\setminus F_{2,v}\r)^-\r).$$
By \eqref{4e11}, following the proof of \eqref{3e14} with some slight modifications, we obtain that,
for any $(u,v)\in B((\widetilde{f},\,\widetilde{g}),d)$ and $y\in V$,
\begin{align*}
\int_K|u|(x)|v|(x^-\ast y)\,d\mu(x)
&>\frac{C_fC_g}4\mu\lf(Y_y\r)\\
&\ge\frac{D^2\mu(F_1)}{72}\vz^{-1}\lf(\frac1{\|\mathbf{1}_{F_1}\|_\Oz}\r)
\phi^{-1}\lf(\frac1{\|\mathbf{1}_{F_2}\|_\Oz}\r)>m.
\end{align*}
This, together with the fact $U\subset V$ gives us $(u,v)\notin E_{U,m}$ and hence finishes the proof
of Theorem \ref{t4}.
\end{proof}

Finally, we give the proof of Theorem \ref{t5}.

\begin{proof}[\textbf{Proof of Theorem \ref{t5}}]
Proceeding as in the previous proof, we will show that, for any $m\in\nn$, the set
$$E_{U,m}:=\lf\{(f,g)\in\cl\times\clp:\ \exists\,x\in U,\
\int_K |f|(y)|g|(y^-\ast x)\,d\mu(y)<m\r\}$$
is nowhere dense.

Note that the compact set $U$ satisfies the condition $(\star\star)$. Then there exist some
$a\in K$ and a compact symmetric neighborhood $V$ of the identity element $e$ in $K$ such that
$\Delta(a)<1$, $U\subset V$ and \eqref{2e9} through \eqref{2e11} are valid. Without loss of generality,
we may assume that
\begin{align}\label{4e12}
\frac{1}{\|\mathbf{1}_{V\ast V}\|_{\Oz}}<\phi(b_{\phi})
\end{align}
since $K$ is not compact and $\phi(b_{\phi})>0$. For any given $D\in(0,\fz)$ and $m\in\nn$, by
\eqref{2e12} and the assumption $\vz(b_{\vz})>0$, we can find some $x_0\in(0,\vz(b_{\vz}))$ such that,
for any $x\in(0,x_0)$,
\begin{align}\label{4e13}
\frac{D^2}{72(M+1)}\frac{\vz^{-1}(x)}{x}
\phi^{-1}\lf(\frac1{C_{(U,V,a)}\|\mathbf{1}_{V\ast V}\|_\Oz}\r)>m,
\end{align}
where $\displaystyle M:=\lim_{s\to\fz}\frac{\Phi_\Oz(s)}s$ and $C_{(U,V,a)}\in[1,\fz)$ is as in \eqref{2e10}.

On another hand, from \eqref{2e11}, the assumption (VI) of $\Oz$ and Lemma \ref{2l1}(iii), it follows that
there exists some $k_0\in\nn$ such that, for any $k\in[k_0,\fz)\cap\nn$,
\begin{align}\label{4e14}
\frac{1}{\|\mathbf{1}_{V\ast\{a^-\}^k}\|_{\Oz}}<x_0\quad {\rm and} \quad
\lf\|\mathbf{1}_{V\ast\{a^-\}^k}\r\|_{\Oz}\le(M+1)\mu\lf(V\ast\{a^-\}^k\r).
\end{align}
Define
$$F_1:=V\ast\{a^-\}^{k_0}\quad {\rm and} \quad
F_2:=\{a\}^{k_0}\ast V\ast V.$$
Then, by \eqref{2e10}, \cite[Lemma 3.3C]{Je75} and \eqref{4e12}, similar to \eqref{3e5}, we have
\begin{align*}
\frac{1}{C_{(U,V,a)}\|\mathbf{1}_{V\ast V}\|_\Oz}\le
\frac{1}{\|\mathbf{1}_{F_2}\|_{\Oz}}\le\frac{1}{\|\mathbf{1}_{V\ast V}\|_{\Oz}}<\phi(b_{\phi}).
\end{align*}
According to this and \eqref{4e14}, the inequality \eqref{4e13} takes the form
\begin{align*}
\frac{D^2\mu(F_1)}{72}\vz^{-1}\lf(\frac{1}{\|\mathbf{1}_{F_1}\|_{\Oz}}\r)
\phi^{-1}\lf(\frac1{\|\mathbf{1}_{F_2}\|_\Oz}\r)>m.
\end{align*}
The remainder of the present proof is parallel to that of Theorem \ref{t4};
we omit the details. The proof of Theorem \ref{t5} is completed.
\end{proof}

\noindent  Jun Liu and Chi Zhang (Corresponding author)

\noindent School of Mathematics, JCAM,
China University of Mining and Technology,
Xuzhou 221116, China

\noindent {\it E-mails}: \texttt{junliu@cumt.edu.cn} (J. Liu);
\texttt{zclqq32@cumt.edu.cn} (C. Zhang)

\medskip

\noindent  Yaqian Lu

\noindent
School of Mathematics and Statistics, Central South University,
Changsha 410075, China

\noindent{\it E-mail:} \texttt{yaqianlu@csu.edu.cn} (Y. Lu)

\end{document}